\documentclass{amsart}
\usepackage[latin1]{inputenc}
\usepackage{amssymb, amsmath}
\usepackage[pdftex]{hyperref}

\numberwithin{equation}{section}

\theoremstyle{plain}
\newtheorem{theorem}{Theorem}[section]

\newtheorem{fact}[theorem]{Fact}
\newtheorem{lemma}[theorem]{Lemma}
\newtheorem{corollary}[theorem]{Corollary}

\theoremstyle{definition}
\newtheorem{definition}[theorem]{Definition}
\newtheorem{remark}[theorem]{Remark}
\newtheorem{hypothesis}[theorem]{Hypothesis}

\def\defn{\begin{definition}}
\def\edefn{\end{definition}}
\def\satz{\begin{theorem}}
\def\esatz{\end{theorem}}
\def\tats{\begin{fact}}
\def\etats{\end{fact}}
\def\kor{\begin{corollary}}
\def\ekor{\end{corollary}}
\def\lmm{\begin{lemma}}
\def\elmm{\end{lemma}}
\def\bem{\begin{remark}}
\def\ebem{\end{remark}}
\def\hyp{\begin{hypothesis}}
\def\ehyp{\end{hypothesis}}
\def\bew{\begin{proof}}
\def\ebew{\end{proof}}

\def\M{\mathfrak M}
\def\N{\mathfrak N}

\def\bdd{\mathrm{bdd}}
\def\acl{\mathrm{acl}}
\def\dcl{\mathrm{dcl}}
\def\aut{\mathrm{Aut}}
\def\Scl{\Sigma\mathrm{cl}}
\def\tp{\mathrm{tp}}

\def\Z{\tilde Z}
\def\St{\mathrm{Stab}}
\def\coker{\mathrm{coker}}

\def\FF{\mathfrak F}

\def\Ind#1#2{#1\setbox0=\hbox{$#1x$}\kern\wd0\hbox to 0pt{\hss$#1\mid$\hss}
\lower.9\ht0\hbox to 0pt{\hss$#1\smile$\hss}\kern\wd0}
\def\Notind#1#2{#1\setbox0=\hbox{$#1x$}\kern\wd0\hbox to
0pt{\mathchardef\nn="0236\hss$#1\nn$\kern1.4\wd0\hss}\hbox to
0pt{\hss$#1\mid$\hss}\lower.9\ht0
\hbox to 0pt{\hss$#1\smile$\hss}\kern\wd0}
\def\ind{\mathop{\mathpalette\Ind{}}}

\begin{document}

\title{Bounded Morphisms}
\date{7 July 2020}

\author{Frank O. Wagner}
\address{Universit\'e de Lyon; Universit\'e Claude Bernard Lyon 1; Institut 
Camille Jordan UMR5208, 43 boulevard du 11
novembre 1918, F--69622 Villeurbanne Cedex, France }
\email{wagner@math.univ-lyon1.fr}
\thanks{Partially supported by ANR-13-BS01-0006 ValCoMo and ANR-17-CE40-0026 Agrume.}
\keywords{bounded automorphism, Pfaffian family, bounded endomorphism}
\subjclass[2010]{03C45, 12H05, 12H10}

\begin{abstract}A bounded automorphism of a field or of a group with
trivial almost-centre is definable. In an expansion of a field by a Pfaffian family $\FF$ of additive endomorphisms such that the algebraic closure of any $\FF$-substructure in the expansion coincides with relative field-theoretic algebraic closure, a bounded ring endomorphism, possibly composed with a power of the Frobenius, is a composition of ring endomorphisms canonically associated to $\FF$.\end{abstract}
\maketitle

\section*{Introduction}
In \cite{La92} Lascar calls an automorphism $\sigma$ of a structure $\M$ {\em bounded} if there is a small set $A$ of parameters such that for all $m\in M$ we have $\sigma(m)\in\acl(A,m)$. He shows that for a strongly minimal set $\M$ the subgroup $\aut_f(\M)$ of strong automorphisms (which fix $\acl^{eq}(\emptyset)$) is simple modulo the normal subgroup of
bounded automorphisms. This result was generalized in \cite{EGT} to structures with an integer dimension compatible with a notion of stationary independence.

More generally, given an invariant family $\Sigma$ of partial types, an automorphism $\sigma$ of a structure $\M$ with simple theory is called {\em$\Sigma$-bounded} if there is a small set $A$ of parameters such that for all $m\in M$ we have $\sigma(m)\in\Scl(A,m)$, where the $\Sigma$-closure $\Scl(B)$ of $B$ is the set of all elements $b$ such that $\tp(b/B)$ is $\Sigma$-analysable. If $\Sigma$ is the class of the algebraic types, then the $\Sigma$-closure is the algebraic closure, and $\Sigma$-bounded coincides with bounded. In \cite{BHMP} Blossier, Hardouin and Mart\'\i n Pizarro show that for a field with
of operators, any $\Sigma$-bounded automorphism is definable and obtained by composition of a power of the Frobenius with automorphisms among the operators or their inverses, where $\Sigma$ is taken as the family of all partial types co-foreign to the generic types of the field. In particular, this yields a uniform proof of results of Lascar \cite{La92}, Ziegler \cite{Zi91} and Konnerth \cite{Ko02}.

Of course, the definition of a bounded automorphism is meaningful for any permutation $\sigma$ of a set $X$; it need not be induced by an automorphism of the ambient structure, nor preserve all the induced structure on $X$. Similarly, the definition of $\Sigma$-bounded is valid for any hyperdefinable set $X$ in a simple theory; moreover we can freely choose $\Sigma$, but we should take $\Sigma$ small with respect to $X$: for instance, any permutation of $X$ is $X$-bounded. As {\em a priori} a $\Sigma$-bounded automorphism of a hyperdefinable set $X^\M$ in a structure $\M$ need not extend to a $\Sigma$-bounded automorphism of $X^\N$ in an elementary extension $\N$ (and even less so with the same small set $A$ of parameters), one should suppose that $\M$ is at least $|T(A)|^+$-saturated, in order not to have to change the model one is working in.

In this paper we shall generalize the main theorem of \cite{BHMP} in two directions. On the one hand, we consider group-theoretic $\Sigma$-bounded automorphisms of a hyperdefinable group $G$ in a simple theory, where $\Sigma$ is still the class of partial types (with parameters) co-foreign to the generic types of $G$. We show that if any non-trivial element has a centralizer of unbounded index (that is, if the almost centre $\Z(G)$ of $G$ is trivial), then any $\Sigma$-bounded automorphism is definable.

On the other hand, we will look at bounded {\em endo}morphisms in an expansion of a field by a family $\FF=\{f_i:i<\alpha\}$ of additive endomorphisms, such that\begin{enumerate}
\item The family $\FF$ is {\em Pfaffian}, that is $f_i(xy)$ is a linear combination of products $f_j(x)f_k(y)$ with $j, k\le i$, for all $i<\alpha$.
\item The algebraic closure of an $\FF$-substructure $A$ in the sense of the expansion is equal to the relative field-theoretic algebraic closure of $A$.\end{enumerate}
These conditions are satisfied in particular by\begin{itemize}
\item differentially closed fields in characteristic zero with several commuting derivations \cite{MG00};
\item generic difference fields in any characteristic \cite{CH99};
\item separably closed fields of finite imperfection degree with a $p$-base named \cite{De88, Wo79};
\item differential generic difference fields in characteristic zero, where the automorphism commutes with the derivation \cite{BM07, BM11};
\item fields with free operators in characteristic zero \cite{MS14}.\end{itemize}
(For the Pfaffian condition, use \cite[Proposition 1.4]{BHMP}.)
We shall show that any bounded endomorphism of $K$ is the composition of endomorphisms canonically associated with functions in $\FF$, which moreover are linear combinations, followed by a power of the inverse of Frobenius. In the cases mentioned above, the endomorphisms (except for the Frobenius for the separably closed fields) will be surjective, and the associated automorphisms will be the automorphisms among the functions in $\FF$, and possibly the Frobenius.

\section{Preliminaries}
In this section we will recall some less known concepts and results for simple groups and theories. For more details, we refer to \cite{Wa00}. As usual, we will work in a monster model of a simple theory, and we fix an $\emptyset$-invariant family $\Sigma$ of
partial types with parameters.
\defn The {\em $\Sigma$-closure\/} $\Scl(A)$ of $A$ is the collection of all hyperimaginaries $a$ such that $\tp(a/A)$ is $\Sigma$-analysable.\edefn
We always have $\bdd(A)\subseteq\Scl(A)$, with equality if and only if $\Sigma$ only contains algebraic types. In general one will choose a set $\Sigma$ of partial types that are small in the sense that the ambient model is not contained in $\Scl(\emptyset)$.
\tats\cite[Lemma 3.5.3]{Wa00} The following are equivalent:\begin{enumerate}
\item $\tp(a/A)$ is foreign to $\Sigma$.
\item $a\ind_A\Scl(A)$.
\item $a\ind_A\dcl(aA)\cap\Scl(A)$.
\item $\dcl(aA)\cap\Scl(A)\subseteq\bdd(A)$.\end{enumerate}\etats
Except in case of equality with the bounded closure, the $\Sigma$-closure has the size of the monster model. The equivalence $(2)\Leftrightarrow(3)$ allows us to use only a small part.
\tats\cite[Lemma 3.5.5]{Wa00}\label{Sindep} If $A\ind_BC$ then $\Scl(A)\ind_{\Scl(B)}\Scl(C)$.
More precisely, for all $A_0\subset\Scl(A)$ we have
$A_0\ind_{B_0}\Scl(C)$, where $B_0=\dcl(A_0B)\cap\Scl(B)$.\etats
More properties of the $\Sigma$-closure can be found in \cite{Wa04, PW13}, as well as a refinement into different levels.

Now let $G$ be a group hyperdefinable over $\emptyset$, and recall that the ambient theory is simple.
\defn\begin{enumerate}\item Two hyperdefinable subgroups $H$ and $K$ are
{\em commensurable} if their intersection is of bounded index in either one.
\item A hyperdefinable subgroup $H$ is {\em locally connected} if it
is equal to any commensurable group-theoretical or automorphic conjugate.\end{enumerate}\edefn
\tats\cite[Corollary 4.5.16, Lemmas 4.5.18 and 4.5.19]{Wa00}\label{normalizer}\begin{enumerate}
\item A locally connected hyperdefinable subgroup of $G$ has a canonical parameter.
\item The normalizer of a locally connected hyperdefinable group is again hyperdefinabe.
\item Any hyperdefinable subgroup $H$ of $G$ has a locally connected component $H^{loc}$, which is the smallest hyperdefinable locally connected subgoup commensurable with $H$. Its canonical parameter is definable over the parameters used to hyperdefine $H$.\end{enumerate}\etats
We also recall the notion of the almost centre.
\defn The almost centre of a hyperdefinable group $G$, denoted $\Z(G)$, is the characteristic subgroup of elements whose centralizer is of bounded index in $G$.\edefn
In a simple theory, the almost centre of $G$ is hyperdefinable over the same parameters as $G$ itself \cite[Proposition 4.4.10]{Wa00}.

Finally, we need Ziegler's lemma (\cite[Theorem 1]{Zi06} for the stable abelian case, and \cite[Lemma 1.2 and Remark 1.3]{BMPW} for the general case).
\tats\label{Stabs} Let $g$ and $g'$ in $G$ be such that $g$, $g'$ and $gg'$ are pairwise independent. Then the left stabilizers of $g$ and $gg'$ are equal, and commensurable with $g\St(g')g^{- 1}$, and all three are $\emptyset$-connected. Moreover, $g$ is generic in the right coset $\St(g)g$, which is $\bdd(\emptyset)$-hyperdefinable, and similarly for $g'$ and $gg'$.\etats

\section{$\Sigma$-bounded group automorphisms}
In this section $G$ will still be a group hyperdefinable over $\emptyset$ in a simple theory; we shall assume that the ambient model is sufficiently saturated. We let $\Sigma$ be the class of all partial types
co-foreign to generic types of $G$. (Since generic types have non-forking extensions that are translates of one another, co-foreign to one is equivalent to co-foreign to all.) In particular, if $G$ has a regular generic type, $\Sigma$ contains of all non-generic partial types.
\satz If $\Z(G)$ is trivial, then any $\Sigma$-bounded automorphism $\sigma$ of $G$ is hyperdefinable.\esatz
\begin{proof} Let $A$ be a set such that $\sigma(g)\in\Scl(A,g)$ for all $g$ in $G$, and recall that we shall need the ambient model to be $|T(A)|^+$-saturated. Consider first $g\in G$ generic over $A$, and $g'\in G$ generic over $A,g,\sigma^{-1}(g)$. Put $A'=\Scl(A)\cap\bdd(A,g',\sigma(g'))$. According to Fact \ref{Sindep} we have $g',\sigma(g')\ind_{A'}\sigma^{-1}(g),g$.
Since generics of $G$ are foreign to $\Sigma$, the elements $g$ and $g'$ remain generic over $A'$. So $g_0=g'\sigma^{-1}(g)$ and $\sigma(g_0)=\sigma(g')g$ are both generic over $A'$, whence over $A$.

Now choose $g_2$ generic in $G$ over $A,g_0$, and put $g_1=g_0^{-1}g_2$.
Then $g_0$, $g_1$ and $g_2 = g_0g_1$ are pairwise independent over $A$.
By Fact \ref{Sindep} their $\Sigma$-closures $\Scl(A, g_0)$, $\Scl(A, g_1)$ and $\Scl(A, g_2)$ are pairwise independent over $\Scl(A)$. In particular, $(g_0,\sigma(g_0))$, $(g_1,\sigma(g_1))$ and $(g_2,\sigma(g_2))$ are pairwise independent over $\Scl(A)$, whence by Fact \ref{Sindep} also over
$$B =\bdd(B)=\Scl(A)\cap\bdd(A, g_0,\sigma(g_0), g_1,\sigma(g_1), g_2,\sigma(g_2)).$$
Moreover, since the generic types of $G$ are foreign to $\Sigma$, the elements $g_0$, $\sigma(g_0)$, $g_1$ and $g_2$
each remain generic over $B$.

For $i\le2$ let $S_i =\St(g_i,\sigma(g_i)/B)$ be the (left) stabilizer in $G\times G$. So $S_i$ is connected over $B$ by Fact
\ref{Stabs}, the right coset $S_i\cdot(g_i,\sigma(g_i))$ is hyperdefinable over $B$, and the pair $(g_i,\sigma(g_i))$ is generic in that coset. Furthermore, $S_0=S_2$; let $S = S_0^{loc}$ be its locally connected component. 

By \cite[Lemma 4.5.5]{Wa00}, since 
$$g_i,\sigma(g_i)\ind_{\bdd(A, g_i,\sigma(g_i))\cap\Scl(A)}B,$$
$S_i$ is commensurable with
$\St(g_i,\sigma(g_i)/\bdd (A, g_i,\sigma(g_i))\cap\Scl(A))$ for $i = 0,2$, whose locally connected component is also $S$. It follows that $S$ and the coset $S\cdot (g_i,\sigma(g_i))$ are hyperdefinable over $\bdd(A,g_i,\sigma(g_i))\cap\Scl(A)$ for $i = 0,2$.
Now $g_2$ was any generic over $A, g_0$. So $S$ is hyperdefinable over $\Scl(A)\cap\bdd(A, a,\sigma(a))$ for any $a$ generic over $A,g_0$, and hence over
$$\bar A =\Scl(A)\cap\bigcap_{a\text{ generic over $A,g_0$}}\bdd (A, a,\sigma(a)).$$
Note that since $g_0$ and $\sigma(g_0)$ are generic over $B$, the
projections of $S_0$, and therefore $S$, to the first and to the second coordinate are generic in $G$, and thus subgroups of bounded index.

Let $(g,g^*)$ be generic in $S$ over $B$. As $(g_0,\sigma(g_0)$ is also generic in $S$ over $B$, there is a generic $(x,y)$ of $S$ over $B,g,g^*$ such that $(gx,g^*y)\equiv_B(g_0,\sigma(g_0))$. Hence $g^*y\in\Scl(A,g,x)$, and $g^*\in\Scl(A,g,x,y)$. But $g^*\ind_{B,g}(x,y)$, whence $$g^*\ind_{\Scl(B,g)}\Scl(B,g,x,y),$$
and $g^*\in\Scl(B,g)$.
 
Put $\coker(S) =\{g\in G: (1, g)\in S\}$ and consider $h\in\coker(S)$. If
$(g, g')\in S$ is generic over $B, h$, then $(g, g'h)$ is still generic over $B, h$, and $h\in\Scl(B, g)$. Since $h\ind_Bg$ we have $h\in\Scl(B) =\Scl(A)$, and $\coker(S)\subseteq\Scl(A)$.

Now let $g''\in G$ be arbitrary, and $g\in G$ generic over $\bar A,g_0,g'',\sigma(g'')$. We put 
$g'= g^{-g''}$ and $B' =\Scl(A)\cap\bdd(A, g, g',\sigma(g),\sigma(g'))$. So $g$ and $g'$ are generic over $B',g_0$, and
$$(g,\sigma(g))\in (S\cdot (g',\sigma(g'))^{(g'',\sigma(g''))}\cap (S\cdot (g,\sigma(g)).$$
As the two cosets $S\cdot (g,\sigma(g))$ and $S\cdot (g',\sigma(g'))$ are hyperdefinable over $B'$,
$$g,\sigma(g)\ind_{B'}g'',\sigma(g'')$$
by Fact \ref{Sindep},
and $(g,\sigma(g))$ is generic in $S\cdot (g,\sigma(g))$ over $B'$,
the two groups $S$ and $S^{(g'',\sigma(g''))}$ are commensurable, and therefore equal by local connectivity. 
It follows that $S$ is normalized by the subgroup $\{(g,\sigma(g)):g\in G\}\le G\times G$. In particular, if $h\in\coker(S)$, then $h^{\sigma(G)}\subseteq\coker(S)\subset\Scl(A)$. 

Now $\sigma$ is an automorphism, whence surjective, and $h^G\subseteq\Scl(A)$. So for $g\in G$ generic over $A, h$ we have $g\ind_{A, h}h^g$,
since $\tp(g/A, h)$ is foreign to $\Sigma$. So $h^g\in\bdd(a,h)$, and there is a generic $g'$ independent of $g$ over $A,h$ with $h^g=h^{g'}$. Thus $g'g^{-1}$ is generic over $A,h$ and centralising $h$, whence
$C_G (h)$ has bounded index in $G$ and $h\in\Z(G) =\{1\}$. It follows that
$\coker(S)$ is trivial.

By Fact \ref{normalizer} the normalizer $N_G(S)$ is hyperdefinable.  Since $(g,\sigma(g))$ normalizes $S$ for all $g\in G$, the projection to the first coordinate of $N_G(S)$ is the whole of $G$. But for $(h,h')\in N_G(S)$ and $(g, g')\in S$ with $g'\in G$ 
generic over $h,\sigma(h),h'$ we have
$$(g^h,g'^{h'})\in S\qquad\mbox{and}\qquad(g^h,g'^{\sigma(h)})\in S^{(h,\sigma(h))}=S.$$
Since $\coker(S)$ is trivial, we get $g'^{h'}= g'^{\sigma(h)}$ and $\sigma(h)h'^{-1}\in C_G (g')$. Thus $g'\in C_G(\sigma(h)h'^{-1})$. But $g'$ is generic over $\sigma(h),h'$, so $\sigma(h)h'^{- 1}\in\Z (G) =\{1\}$ and
$h'=\sigma(h)$. Thus $N_G(S)$ hyperdefines $\sigma$.
\end{proof}

\kor\label{bound} Let $\sigma$ be a $\Sigma$-bounded automorphism of a field $K$
hyperdefinable in a simple theory, where $\Sigma$ is the class of partial types co-foreign to the generic types of $K$. Then $\sigma$ is hyperdefinable.\ekor
\begin{proof} Let $G = K^+\rtimes K^\times$. So $G$ is hyperdefinable, and $\Sigma$ is also the set of partial types co-foreign to generics of $G$ (which are pairs of two independent generics of $K$). Furthermore, $\Z(G) =\{1\}$, and $\sigma$ induces an automorphism $\tau: (a, b)\mapsto (\sigma(a),\sigma(b))$ of $G$. According to the previous theorem
$\tau$ is hyperdefinable, and so is $\sigma$.
\end{proof}

\section{Bounded endomorphisms of a field with operators}
In this section we shall consider a hyperdefinable field $K$ in a simple theory, together with a well-ordered family $\FF =\{f_i: i <\alpha\}$ of additive endomorphisms. We assume that\begin{enumerate}
\item The family $\FF$ is $K$-free: if $\sum_i\lambda_if_i = 0$ where $\lambda_i\in K$ is zero for almost all $i<\alpha$, then $\lambda_i = 0$ for all $i <\alpha$.
\item The family is {\em$K$-Pfaffian}, i.e.\ for all $i <\alpha$ there are $a^i_{j, k}\in K$ with $j,k\le i$, almost all $0$, such that
$$f_i(xy) =\sum_{j,k\le i}a^i_{j, k}f_j(x) f_k(y).$$\end{enumerate}
\bem If we start with a $K$-Pfaffian family $\FF$ that is not $K$-free, we can extract a $K$-free and $K$-Pfaffian subfamily which generates the whole of $\FF$. In particular, in positive characteristic we can assume that the Frobenius is part of $\FF$.\ebem
\lmm\label{unique} If $\sum_{\bar\iota}a_{\bar\iota}f_{i_0}(x_0)\cdots f_{i_n}(x_n) =\sum_{\bar\iota}b_{\bar\iota}f_{i_0}(x_0)\cdots f_{i_n}(x_n)$ where $a_{\bar\iota}, b_{\bar\iota}\in K$ are almost all zero, then $a_{\bar\iota}= b_{\bar\iota}$ for all $\bar\iota = (i_0,\ldots, i_n)$.\elmm
\bew By induction on $n$. If $n =$0, this is just $K$-freeness of $\FF$. We therefore assume the assertion true for $n-1$, and
$$\sum_{\bar\iota}a_{\bar\iota}f_{i_0}(x_0)\cdots f_{i_n}(x_n) =\sum_{\bar\iota}b_{\bar\iota}f_{i_0}(x_0)\cdots f_{i_n}(x_n).$$
So
$$0 =\sum_{i_n}\big (\sum_{i_0,\ldots, i_{n-1}}(a_{\bar\iota}-b_{\bar\iota}) f_{i_0}(x_0)\cdots f_{i_{n-1}}(x_{n-1})\big) f_{i_n}(x_n).$$
If there is some $\bar k$ with $a_{\bar k}\not = b_{\bar k}$, then by induction hypothesis
$$g_{i_n}(x_0,\ldots, x_{n-1}) =\sum_{i_0,\ldots, i_{n-1}}(a_{\bar\iota}-b_{\bar\iota}) f_{i_0}(x_0)\cdots f_{i_{n-1}}(x_{n-1})$$
is not identically $0$ for $i_n=k_n$, and there is $\bar\alpha\in K^n$ with $g_{k_n}(\bar\alpha)\not = 0$. But then $\sum_{i_n}g_{i_n} (\bar\alpha) f_{i_n}(x_n) = 0$, which again contradicts $K$-freeness of $\FF$.\ebew
\kor The coefficients $a^i_{j, k}$ are uniquely determined.\ekor
\bew Immediate.\ebew
\bem $a^i_{j,k}=a^i_{k,j}$ for all $j,k\le i$, and $\sum_{j\le i}a^i_{j,k}f_j(1)=\delta_{ki}$ for all $k\le i$.\ebem
\bew The first equation follows from $f_i(xy)=f_i(yx)$, and the second from $f_i(1y)=f_i(y)$.\ebew
\kor\label{product} For all $i <\alpha$ the function
$$\sigma_i(x) =\sum_{j\le i}a^i_{j, i}f_j (x)$$ is a ring endomorphism of $K$. Moreover, we can assume $\sigma_i\in\FF$.\ekor
\bew We have
$$\begin{aligned}f_i(xy)&=\sum_{j,k\le i}a^i_{j,k}f_j(x)f_k(y)=\sum_{j\le i}a^i_{j,i}f_j(x)f_i(y)+\sum_{j\le i,\,k<i}a^i_{j,k}f_j(x)f_k(y)\\
&=\sigma_i(x)f_i(y)+R(x,y),\end{aligned}$$
where $R(x,y)=\sum_{j\le i,\,k<i}a^i_{j,k}f_j(x)f_k(y)$.
Thus
$$\begin{aligned}f_i(xx'y)&=\sigma_i(xx')f_i(y)+R(xx',y)\\
&=\sigma_i(xx')f_i(y)+\sum_{j\le i,\,k<i}a^i_{j,k}\big(\sum_{\ell,m\le j}a^j_{\ell,m}f_\ell(x)f_m(x')\big)f_k(y)\\
&=\sigma_i(xx')f_i(y)+\sum_{j\le i,\,k<i}\ \sum_{\ell,m\le j}a^i_{j,k}a^j_{\ell,m}f_\ell(x)f_m(x')f_k(y)\qquad\end{aligned}$$
and
$$\begin{aligned}
f_i(xx'y)&=\sigma_i(x)f_i(x'y)+R(x,x'y)\\
&=\sigma_i(x)\sigma_i(x')f_i(y)+\sigma_i(x)R(x',y)+R(x,x'y)\\
&=\sigma_i(x)\sigma_i(x')f_i(y)+\sum_{j\le i}a^i_{j, i}f_j (x)\ \sum_{\ell\le i,\,m<i}a^i_{\ell,m}f_\ell(x')f_m(y)\\
&\qquad\qquad\qquad\qquad+\sum_{j\le i,\,k<i}a^i_{j,k}f_j(x)\sum_{\ell,m\le k}a^k_{\ell,m}f_\ell(x')f_m(y)\\
&=\sigma_i(x)\sigma_i(x')f_i(y)+\sum_{j\le i}\ \sum_{\ell\le i,\,m<i}a^i_{j, i}a^i_{\ell,m}f_j (x)f_\ell(x')f_m(y)\\
&\qquad\qquad\qquad\qquad+\sum_{j\le i,\,k<i}\ \sum_{\ell,m\le k}a^i_{j,k}a^k_{\ell,m}f_j(x)f_\ell(x')f_m(y).\end{aligned}$$
Comparing the coefficient of $f_i(y)$, Lemma \ref{unique} yields
$$\sigma_i(xx')=\sigma_i(x)\sigma_i(x')$$
and $\sigma_i$ is a ring endomorphism of $K$. Note that $\sigma_i(1)=1$, so $\sigma_i$ is non-zero.

The $K$-linear combinations of functions in $\FF$ form a $K$-vector space with basis $\FF$. Since every family of (non-zero) ring endomorphisms of $K$ is linearly independent, we can start another basis with $(\sigma_i:i<\alpha)$ (eliminating any repetitions which might occur), followed by the functions of $\FF$, eliminating those which are already in the span of the previous ones together with $(\sigma_i:i<\alpha)$.
It is clear that this family is still Pfaffian, and $K$-free.\ebew
\defn\label{multiplicative} We call $\sigma_i(x)=\sum_{j\le i}a^i_{j,i}f_j(x)$ the {\em endomorphism associated to~$f_i$}.\edefn
Note that the endomorphism associated with a ring endomorphism $\sigma$ is $\sigma$ itself, and the endomorphism associated with a derivation is the identity.
\bem\cite[Corollary 1.5]{BHMP}\label{mult} If we well-order compositions of additive endomorphisms in $\FF$ first by the length of the composition and then by lexicographic order of the indices of the functions of $\FF$ used, the endomorphism associated to such a composition will be the composition of the associated endomorphisms of the functions in $\FF$. We shall denote this well-ordered sequence of additive endomorphisms by $(f_\theta)_\theta$. The proof of Corollary \ref{product} yields that $f_\theta(ax)=\sigma_\theta(a)f_\theta(x)+\sum_{\theta'<\theta}a^\theta_{\theta'}f_{\theta'}(x)$.\ebem

The following lemma will not be used for the proof of the main theorem of this section.
\lmm Any finite subfamily $\{f_i: i\in I_0\}$ of $\FF$ is contained in a finite Pfaffian subfamily $\{f_i: i\in I\}$, with $\max I =\max I_0$.\elmm
\bew We put $I_{- 1}=\emptyset$, and for $n> 0$ we build inductively a string of finite families $I_0\subseteq I_1\subseteq I_2\subseteq\cdots$ by posing
$$I_{n + 1}= I_n\cup\{j <\alpha:\exists\, i\in I_n\,\exists k\, a^i_{j, k}\not = 0\}.$$
Since the sequence $\max (I_n\setminus I_{n-1})$ is strictly decreasing, the sequence of the $I_n$ becomes stationary, and $I =\bigcup_{n <\omega }I_n$ gives a Pfaffian family $\{f_i: i\in I\}$.\ebew

\satz\label{expansion} Let $K$ be a field, and $\FF$ a Pfaffian family of additive endomorphisms of $K$. For $A\subseteq K$ let $\langle A\rangle$ be the closure of $A$ under all functions $f\in\FF$. We assume that $\acl(A) =\langle A\rangle^{alg}\cap K$, where $\acl$ is the algebraic closure in the sense of the structure $(K, +,\cdot, f: f\in\FF)$, and $A^{alg}$ is the field-theoretic algebraic closure. Then any definable endomorphism $\sigma$ of the field $K$ satisfies an identity of the form
$$\mbox{Frob}^j\circ\sigma=\sigma_i,$$
where $\sigma_i$ is a composition of associated endomorphisms of $\FF$, and the Frobenius (where $\mbox{Frob}^0$ is the identity, and $j=0$ in characteristic $0$).\esatz
\bew We can assume that $\FF$ is $K$-free, and that in positive characteristic the Frobenius is in $\FF$. Let $B =\acl (B)$ be parameters over which $\sigma$ is definable.

For all $a\in K$ the image $\sigma(a)$ is in 
$$\acl (Ba) =\langle B a\rangle^{alg}\cap K = (B\cup\langle a\rangle)^{alg}\cap K.$$
By compactness, there is a finite number of polynomials $(P_i(\bar x_i, y): i <\ell)$ over $B$, and for all $i <\ell$ a sequence $(f_{\theta_{i, j}}: j <|\bar x_i |)$ of compositions of functions in $\FF$ such that for all $a\in K$ there is $i <\ell$ with $P_i(\bar f_{\theta_{i, j}}(a), y)$ non-trivial in $y$, and satisfied by $\sigma(a)$.

Let $\bar K = K^{alg}$ be the algebraic closure of $K$, and $(f_{\theta_i}: i <n)$ the sequence of functions $f_{\theta_{i, j}}$, with $f_0 =\,$id. We consider the following additive subgroup of $\bar K^{n + 1}$:
$$G =\{((f_{\theta_i}(a): i <n),\sigma(a)): a\in K\}.$$
Let $\Gamma$ be its Zariski closure, and $\Gamma^0$ the connected component of $\Gamma$. Then $\Gamma^0$ is the Zariski closure of $G^0 =\Gamma^0\cap G$, and the latter is of finite index in $G$.

Since every element $(\bar x,y)$ of $G$ satisfies an equation $P_i(\bar x, y) = 0$ (with a suitable enumeration of $\bar x$), non-trivial in $y$, for some $i <\ell$, the disjunction of these equations is generically satisfied in $\Gamma^0$ and non-trivial in $y$ by \cite[Remark 27]{Wa90}.

However, the only proper additive algebraic subgroups of $\bar K^{n + 1}$ are given by additive polynomials. So there is a polynomial $P\in\bar K[X_0,\ldots, X_n]$ of the form
$$P(\bar x, x_n) =\sum_{i\le n,\, j <\omega}\lambda_{i, j}x_i^{p^j}$$
such that $P(\bar x, x_n) = 0$ is satisfied on $\Gamma^0$, and such that $P(\bar 0, x_n)$ is non-trivial. (In characteristic zero we have $j=0$, so $P$ is linear in each $x_i$.)

Since $G$ and $G^0$ are Gal$(K^{sep}/K)$-invariant, so are $\Gamma$ and $\Gamma^0$, and we can assume that $P$ has coefficients in the purely inseparable closure of $K$. By composing with a power of Frobenius, we can even suppose that $P\in K[X_0,\ldots, X_n]$.

The function $x\mapsto P ((f_{\theta_i}(x))_{i <n},\sigma(x))$ is an additive endomorphism of $K$ whose image $F$ is a finite additive subgroup. In characteristic zero $F$ is trivial; in positive characteristic there is an additive nontrivial polynomial $Q\in K[X]$ which vanishes on $F$. So
$$Q(P ((f_{\theta_i}(x))_{i <n},\sigma(x))) = 0$$ for all $x\in K$. Since $\FF$ contains Frobenius in positive characteristic, we obtain a non-trivial equation of the form
$$\sum_{i <n}\lambda_i f_{\theta_i}(x) =\sum_{j <\omega}\mu_j\sigma(x)^{p^j}$$
satisfied on $K$ (possibly with $n$ and $\theta_i$ different). We choose such an equation with
$$\theta_m =\max\{\theta_i:\lambda_i\not = 0\}$$
minimal. Moreover, we can assume that there is only one non-zero $\mu_j$, since for $j_0$ maximal with $\mu_{j_0}\not = 0$ and $a\in K$ transcendental,
$$\sum_{i <n}\lambda_i f_{\theta_i}(ax) -\sigma(a)^{p^{j_0}}\sum_{i <n}\lambda_i f_{\theta_i}(x) =\sum_{j <j_0}\mu_j [\sigma(a)^{p^j}-\sigma(a)^{p^{j_0}}]\,\sigma(x)^{p^j}$$
allows us to reduce the number of non-trivial $\mu_j$ while not increasing $\theta_m$; note that since $a$ and $\sigma(a)$ are transcendental, the equation remains non-trivial if there are at least two non-zero $\mu_j$.

We can now take $\mu = 1$, and obtain $\sum_{i <n}\lambda_i f_{\theta_i}(x) =\sigma(x)^{p^j}$ for all $x\in K$.
For all $a\in K$ we have from Remark \ref{mult} that
$$\begin{aligned}\sum_{i <n}\lambda_i f_{\theta_i}(ax) & =\lambda_m\sigma_{\theta_m}(a) f_{\theta_m}(x) +\sum_{\theta <\theta_m}\alpha_\theta f_\theta (x)\\
& =\sigma(ax)^{p^j}=\sigma(a)^{p^j}\sigma(x)^{p^j}\end{aligned}$$
for some coefficients $\alpha_\theta\in K$ that depend on $a$. If $\sigma_{\theta_m}(a)\not =\sigma(a)^{p^j}$, then
$$\begin{aligned}{}[\sigma(a)^{p^j}-\sigma_{\theta_m}(a)]\,\sigma(x)^{p^j}&=
\lambda_m\sigma_{\theta_m}(a) f_{\theta_m}(x) +\sum_{\theta<\theta_m}\alpha_\theta f_\theta (x) -\sigma_{\theta_m}(a)\sum_{i <n}\lambda_i f_{\theta_i}(x)\\
& = \sum_{\theta<\theta_m}\alpha_\theta f_\theta (x) -\sum_{i <n,\,i\not=m}
\sigma_{\theta_m}(a)\lambda_i f_{\theta_i}(x)\end{aligned}$$
gives a non-trivial equation with a smaller maximal $\theta$, a contradiction.
The theorem follows.\ebew

This proof is inspired by that of Blossier, Hardouin and Martin Pizarro \cite[Th\'eor\`eme 3.1]{BHMP}, but it uses the generic properties of the definable envelope (which here is equal to the Zariski closure) of \cite{Wa90} instead of the stabilizer. We need for our Theorem \ref{expansion} only part of hypothesis 4 of \cite{BHMP} (the characterization of the algebraic closure) and we make no assumptions about the simplicity of the expansion; on the other hand, we assume from the beginning the Pfaffian property (\cite[Proposition 1.4]{BHMP}, which follows from their hypotheses 1 and 2). Note that we are not considering the inverses of automorphisms in $\FF$: we can always add them to $\FF$ and preserve the Pfaffian property; and in the applications this will be necessary to obtain the characterization of the algebraic closure.

We can now obtain Th\'eor\`eme 3.1 of \cite{BHMP} (note that {\em bounded} there means {\em$\Sigma$-bounded} in our sense). Firstly, by Corollary \ref{bound} any $\Sigma$-bounded automorphism of a simple field must be definable. Secondly, under hypotheses 1 and 2 of \cite[Th\'eor\`eme 3.1]{BHMP}, the family of automorphisms can be replaced by a Pfaffian family by \cite[Proposition 1.4]{BHMP}. Hypothesis 5 of \cite[Th\'eor\`eme 3.1]{BHMP} tells us that $\acl(A) =\langle A\rangle^{alg}\cap K$. Now by Theorem~\ref{expansion} there is an identity $\mbox{Frob}^j\circ\sigma=\sigma_i$,
where $\sigma_i$ is a composition of associated endomorphisms of $\FF$, and the Frobenius. But $\sigma_i$ is surjective by hypothesis 2, so either $j=0$ or the Frobenius is invertible. In either case, $\sigma$ is a composition of associated endomorphisms, the Frobenius, and its inverse if applicable.

\end{document}